\documentclass[12pt]{article}
\usepackage{amsmath}
\usepackage{graphicx}
\usepackage{amsfonts}
\usepackage{amssymb}
\newcommand{\cc}{{\cal C}}

\newcommand{\R}{{\mathbb R}}

\newtheorem{theorem}{Theorem}[section]

\newtheorem{lemma}[theorem]{Lemma}

\begin{document}
\thispagestyle{empty}
 \begin{center}
 {\Large\bf Iterated logarithm law for anticipating stochastic differential equations}

\vspace{0,2cm}

by

\vspace{0,2cm}

{\bf David M\'arquez-Carreras}\footnote{Corresponding
author.\hfill}$^{,2}$ and {\bf Carles Rovira}\footnote{Partially
supported by DGES grant BFM2003-01345.\hfill}

\vspace{0,1cm}

Facultat de Matem\`atiques, Universitat de Barcelona,
\\
Gran Via 585, 08007-Barcelona, Spain
\\
e-mail: davidmarquez@ub.edu, carles.rovira@ub.edu

 \end{center}

\vspace{0.3cm}

\begin{abstract}
We prove a functional law of iterated  logarithm for the following
kind of anticipating stochastic differential equations
$$\xi^u_t=X_0^u+\frac{1}{\sqrt{\log\log u}}\sum_{j=1}^k \int_0^{t}
A_j^u(\xi^u_s)\circ dW_{s}^j+ \int_0^{t} A_0^u(\xi^u_s)ds,$$ where
$u>e$, $W=\{(W_t^1,\dots,W_t^k), 0\le t\le 1\}$ is a standard
$k$-dimensional Wiener process, $A_0^u,A_1^u,\dots,
A_k^u:\mathbb{R}^d\longrightarrow \mathbb{R}^d$ are functions of class
$\mathcal{C}^2$ with bounded partial derivatives up to order $2$,
$X_0^u$ is a random vector not necessarily adapted and the first
integral is a
generalized Stratonovich integral . %In this paper, more
%specifically, we check that the family $\{\xi^u_\cdot,\ u>e\}$ is
%relatively compact and that the almost surely limit of this family
%when $u$ goes to infinity exists, identifying this limit.
\end{abstract}

\vspace{0,5cm}

\noindent {\bf Running head:} ILL for anticipating SDE



\newpage

\section{Introduction}

Consider the
 Stratonovich differential equation on $\mathbb{R}^d$
\begin{equation}
X_t=X_0+\sum_{j=1}^k \int_0^t A_j(X_s)\circ dW_s^j+\int_0^t
A_0(X_s)ds, \label{esdesen}
\end{equation}
where $W=\{(W_t^1,\dots,W_t^k), 0\le t\le 1\}$ is a standard
$k$-dimensional Wiener process, $A_0,A_1,\dots,
A_k:\mathbb{R}^d\longrightarrow \mathbb{R}^d$ are functions of class
$\mathcal{C}^2$ with bounded partial derivatives up to order $2$
and $X_0$ is a random vector not necessarily adapted to the
filtration associated with the Wiener process. Here the stochastic
integral term is defined as an anticipating Stratonovich integral
(see, for instance, the paper of Nualart and Pardoux, 1988). Under
some smooth conditions on $X_0$ and the coefficients, Ocone and
Pardoux (1989) prove the existence and uniqueness of solutions for
(\ref{esdesen}).

Millet, Nualart and Sanz-Sol\'{e} (1992) consider, for
$\varepsilon>0$, the following family of perturbed anticipating
stochastic differential equation
\begin{equation}
X_t^\varepsilon=X_0^\varepsilon+\sqrt{\varepsilon}\sum_{j=1}^k
\int_0^t A_j(X_s^\varepsilon)\circ dW_s^j+\int_0^t
A_0(X_s^\varepsilon)ds. \label{esdep}
\end{equation}
They show that a solution of (\ref{esdep})  can be expressed as
the composition of the following adapted flow
\begin{equation}
\varphi_t^\varepsilon(x)=x+\sqrt{\varepsilon}\sum_{j=1}^k \int_0^t
A_j(\varphi_s^\varepsilon(x))\circ dW_s^j+\int_0^t
A_0(\varphi_s^\varepsilon(x))ds, \quad x\in
\mathbb{R}^d\label{eflow}
\end{equation}
and the initial condition, that means
$X_t^\varepsilon=\varphi_t^\varepsilon(X_0^\varepsilon)$. They
also obtain a large deviatons principle (LDP) for the family of
laws of $\{X^\varepsilon\}_{\varepsilon > 0}.$

 It is natural to study the existence of an almost sure functional
 law of iterated logarithm generalizing the Strassen Theorem. This
 problem  has been studied for diffusions by Baldi (1986), for
 parabolic SPDEs by Chenal and Millet (1999) and for
 stochastic Volterra equations by Ait Ouahra and Mellouk (2005).
In this paper, following the ideas presented by Baldi (1986), we
prove a similar result for an anticipating stochastic differential
equation.

The structure of the paper is the following. In Section 2 we
recall some notations and results of Millet, Nualart and Sanz
(1992) about the large deviations principle for anticipating
stochastic differential equations. In Section 3 we present our
equation and we adapt the results of Millet, Nualart and Sanz
(1992) to our framework. Finally, in Section 4, we present our law
of iterated logarithm.

\section{Large deviations principle}\label{}

In order to present a large deviation principle we borrow the
notations of Millet, Nualart and Sanz (1992). For any integer
$m\ge 1$ and $x\in \mathbb{R}^m$, we denote by $H_x^m$ the set of
absolutely continuous functions $f\in
\mathcal{C}([0,1];\mathbb{R}^m)$ with $f_0=x$ and $\int_0^1 |\dot
f_s|^2\ ds<+\infty$. If $x=0$ we write $H^m$ instead of $H_0^m$.
Given $f\in H^k$ we consider the function $g(x)\in H_x^d$, which
is the solution of the differential equation
\begin{equation}\label{edeterm}
g_t(x)=x+\sum_{j=1}^k \int_0^t A_j(g_s(x))\ \dot f^j_s\ ds
+\int_0^t A_0(g_s(x))ds.
\end{equation}
Millet, Nualart and Sanz-Sol\'{e} (1992) prove the following
Theorem:

\newpage

\begin{theorem}\label{tpgd}
Assume that:\begin{itemize} \item[{\rm (h)}] The coefficients
$A_0, A_1,\dots, A_k, B$ and $M=\frac12 \sum_{j=1}^k A_j
\partial A_j$ are of class $\mathcal{C}^2$ with bounded partial
derivatives up to order $2$. \item[{\rm (c)}] There exists $x_0\in
\mathbb{R}^d$ such that for any $\delta>0$
$$\limsup_{\varepsilon\rightarrow 0} \varepsilon \log
P\{|X_0^\varepsilon-x_0|>\delta\}=-\infty.$$\end{itemize} Then,
the family $\{P^\varepsilon,\ \varepsilon>0\}$ of laws of
$\{X^\varepsilon_\cdot=\varphi_\cdot^\varepsilon(X_0^\varepsilon),\
\varepsilon>0\}$ satisfies a large deviation principle with rate
function
\begin{equation}\label{erate}
I(g)=\inf\{\mathcal{I}(f);\ f\in H^k,\
g=F_{x_0}(f)\},\end{equation} where $F_{x_0}(f)$ denotes the
solution of the ordinary differential equation {\rm
(\ref{edeterm})} with initial condition $x=x_0$ and
\begin{equation}\label{erate2}
 \mathcal{I}(f)=\left\{\begin{array}{ll} \frac12 \int_0^1
|\dot
f_s|^2\ ds, & {\rm if} \ f\in H^k,\\
+ \infty, & {\rm otherwise}.\end{array}\right.\end{equation}
\end{theorem}

\section{Structure of our equation}

We denote by $\varphi_t(x)$ the flow $\varphi_t^\varepsilon(x)$ of
(\ref{eflow}) when $\varepsilon=1$. Following the methods
introduced in Millet, Nualart and Sanz (1992), it will be useful
to express now $\varphi_t(x)$ using It\^{o} integral. So, we can
rewrite $\varphi_t(x)$ in the following form
\begin{equation}\label{eflow2}
\varphi_t(x)=x+\sum_{j=1}^k \int_0^t A_j(\varphi_s(x))\
dW_s^j+\int_0^t B(\varphi_s(x))ds,
\end{equation}
with $B=A_0+\frac12 \sum_{j=1}^k A_j \partial A_j$ and where the
stochastic integral term is now defined as an It\^{o} integral.

For $u>e$ we define
\begin{equation*}\phi(u)=\sqrt{u L(u)},\qquad {\rm with}\
L(u)= \log\log u.
\end{equation*}
Let $\mu^u_t(x)=\phi(u)^{-1} \varphi_{u t}( \phi(u) x)$. Using a
change of variable and the scaling property we have that
\begin{equation*}\begin{array}{l}
\displaystyle \mu^u_t \left(\frac{x}{\phi(u)}
\right)=\frac{x}{\phi(u)}+\frac{1}{\phi(u)}\left(\sum_{j=1}^k
\int_0^{u t} A_j(\varphi_s(x))\ dW_s^j+\int_0^{u t}
B(\varphi_s(x))ds\right)\\[5mm]
\displaystyle \qquad
=\frac{x}{\phi(u)}+\frac{1}{\phi(u)}\sum_{j=1}^k \int_0^{t}
A_j(\varphi_{u s}(x))\ dW_{u s}^j+\frac{u}{\phi(u)}
\int_0^{t} B(\varphi_{u s}(x))ds\\[5mm]
\displaystyle \qquad
=\frac{x}{\phi(u)}+\frac{\sqrt{u}}{\phi(u)}\sum_{j=1}^k \int_0^{t}
A_j(\varphi_{u s}(x))\ d \hat W_{s}^j+\frac{u}{\phi(u)} \int_0^{t}
B(\varphi_{u s}(x))ds,
\end{array}\end{equation*}
where $\hat W$ denotes a standard $k$-dimensional Wiener process
that we will also denote by $W$. Then, we can write
\begin{eqnarray*}\label{eoureq}
\mu^u_t\left(\frac{x}{\phi(u)}
\right)&=&\frac{x}{\phi(u)}+\frac{1}{\sqrt{\log\log
u}}\sum_{j=1}^k \int_0^{t} A_j^u
\left(\mu^u_s\left(\frac{x}{\phi(u)}
\right)\right)\circ dW_{s}^j \\
& &  + \int_0^{t} A_0^u \left(\mu^u_s\left(\frac{x}{\phi(u)}
\right)\right)ds,
\end{eqnarray*}
where
\begin{eqnarray*}A_j^u(z)&=& A_j(\phi(u)z),\quad  j=1,\dots,k,\\
A_0^u(z)&=&\frac{u}{\phi(u)}\left[ B(\phi(u) z)-\frac12
\sum_{j=1}^k \sum_{l=1}^d (A_j)^l \partial_l A_j (\phi(u)
z)\right].
\end{eqnarray*}
Consider now the stochastic flow
\begin{equation}\label{flowbo}
\eta^u_t(x)=x+\frac{1}{\sqrt{\log\log u}}\sum_{j=1}^k \int_0^{t}
A_j^u (\eta^u_s(x))\circ dW_{s}^j+ \int_0^{t} A_0^u
(\eta^u_s(x))ds.
\end{equation}
Denote $X_0^u=\phi(u)^{-1} X_0$ and $\xi^u_t \equiv
\eta_t^u(X_0^u)$. Notice that under nice conditions on the
coefficients (see for instance Theorem 3.1 in Millet, Nualart and
Sanz, 1992), we have
\begin{equation}\label{eoureq}
\xi^u_t=X_0^u+\frac{1}{\sqrt{\log\log u}}\sum_{j=1}^k \int_0^{t}
A_j^u(\xi^u_s)\circ dW_{s}^j+ \int_0^{t} A_0^u(\xi^u_s)ds.
\end{equation}
We can now state the following theorem.
\begin{theorem}\label{tpgdnostre}
Assume that: \begin{itemize} \item[{\rm (H)}] The coefficients
$A_0^u, A_1^u,\dots, A_k^u$ and $M^u=\frac12 \sum_{j=1}^k A_j^u
\partial A_j^u$ are of class $\mathcal{C}^2$ with bounded partial
derivatives up to order $2$ and there exist $\tilde A_0, \tilde
A_1$, $\dots,\tilde A_k$ of class $\mathcal{C}^1$ such that
$$\lim_{u\rightarrow + \infty} A_j^u(x)=\tilde A_j(x),\qquad
\lim_{u\rightarrow + \infty} \partial  A_j^u(x)=
\partial \tilde A_j(x),\qquad \forall j=0, 1,\dots,k, $$
uniformly on compact sets on $\mathbb{R}^d$. \item[{\rm (C)}] For
any $\delta>0$,
$$\limsup_{u\rightarrow + \infty} \frac{1}{\log\log u} \log
P\{|X_0^u|>\delta\}=-\infty.\quad  $$\end{itemize} Then, the
family $\{P^u,\ u>e\}$ of laws of $\{\xi^u_\cdot,\ u>e\}$
satisfies a large deviation principle with rate function
\begin{equation}\label{erate}
\tilde I(\tilde g)=\inf\{\mathcal{I}(f);\ f\in H^k,\ \tilde
g=\tilde F_{0}(f)\},\end{equation} where $\mathcal{I}(f)$ is
defined in {\rm (\ref{erate2})} and $\tilde F_{0}(f)$ denotes the
solution of the ordinary differential equation
\begin{equation}\label{edetermnostre}
\tilde g_t=\sum_{j=1}^k \int_0^t \tilde A_j(\tilde g_s)\ \dot
f^j_s\ ds +\int_0^t \tilde A_0(\tilde g_s)ds.
\end{equation}
\end{theorem}
{\bf Proof}: The same proofs as in Millet, Nualart and Sanz (1992)
changing $A_0^u, A_1^u, \dots, A_k^u, M^u$ by $A_1,\dots,A_k$,
$B$, $M$, respectively, still work. The proof is based on an
inequality that we will give in Theorem \ref{desnostre}. \hfill
$\Box$

\medskip

Recall the family $\{P^u,\ u>e\}$ of laws of $\{\xi^u_\cdot,\
u>e\}$ satisfies a large deviation principle with rate function
$\tilde{I}$ defined in (\ref{erate}) if $\tilde{I}$ is lower
semicontinuous; for every $a>0$ the set $\{\tilde g\in
\mathcal{C}([0,1];\mathbb{R}^d);\ \tilde{I}(\tilde g)\le a\}$ is
compact; and for any open set $G$ and any closed $F$ of the space
$\mathcal{C}([0,1];\mathbb{R}^d)$.
\begin{equation}\label{eliminf}
\liminf_{u\rightarrow + \infty} \frac{1}{\log\log u} \log
P\{\xi^u_\cdot\in G\}\ge - \inf\{\tilde{I}(\tilde g), \tilde g\in
G\},\end{equation} and
\begin{equation}\label{elimsup}
\limsup_{u\rightarrow + \infty} \frac{1}{\log\log u} \log
P\{\xi^u_\cdot\in F\}\le - \inf\{\tilde{I}(\tilde g), \tilde g\in
F\}.\end{equation}

\medskip

In the sequel we will denote by $\left\| \cdot \right\|$ the
supremum norm on $\cc([0,1];\R^d)$ and by $\left\| \cdot
\right\|_O$ the supremum norm on $\cc([0,1] \times O;\R^d)$ for
$O\subseteq \R^d$.

\begin{theorem}\label{desnostre}
 Assume {\rm (H)}. Fix $\lambda>0$. Then, for every
 positive reals $R, \tau$ and a compact subset
$O$ of $\R^d$, there exist $u_0>e$ and $\alpha>0$ such that, for
any $u \ge u_0$ and $f\in H^k$ with $\mathcal{I} (f) \le \lambda$,
we have
$$P\left(   \left\| \eta^{u} - \tilde g \right\|_O
>\tau, \left\|\frac{1}{\sqrt{\log\log u}}W - f
\right\| \le \alpha \right) \le \exp \left(-R \log\log u \right),
$$ where,  for $f\in H^k$, $\tilde g_t$ is the solution of the ordinary
differential equation {\rm (\ref{edetermnostre})} and $\eta_t^u$
is the adapted flow defined by {\rm (\ref{flowbo})}.\end{theorem}
{\bf Proof}: Again, the proof follows the computations given in
Millet, Nualart and Sanz (1992) . \hfill $\Box$

\section{The law of iterated logarithm}

The main result of this paper is the following.
\begin{theorem}\label{tlil}
Assuming {\rm (H)} and {\rm (C)}, the family $\{\xi^u_\cdot,
u>e\}$ is relatively compact. Moreover, the a.s. limit set of
$\{\xi^u\}$ when $u$ goes to infinity is
$$\Theta=\{ \tilde g \in \mathcal{C}([0,1];\mathbb{R}^d); \tilde{I}(\tilde
g)\le 1\}.$$\end{theorem}

In order to prove this theorem, following the method presented by
Baldi (1986), we need to check some preliminary lemmas. For the
sake of completeness, we will give the main arguments.

\begin{lemma}\label{l2.3} For every $c>1$ and $\rho>0$, there
exists a.s. $i_0=i_0(\omega)\in \mathbb{N}$ such that for every
$i>i_0$ we have
\begin{equation}\label{el2.3}
d(\xi^{c^i},\Theta):=\inf_{\tilde g\in \Theta} d(\xi^{c^i},\tilde
g)<\rho,\end{equation} where
$$d(\xi^{c^i},\tilde g)=\sup_{0\le t\le
1}|\xi^{c^i}_t-\tilde g_t|.$$ \end{lemma} {\bf Proof}: Let
$\Theta_\rho=\{\tilde g\in \mathcal{C}([0,1];\mathbb{R}^d);
d(\tilde g,\Theta)\ge \rho\}$. We first prove that there exists
$\delta>0$ such that
\begin{equation}\label{el2.3}
\inf_{\tilde g\in \Theta_\rho} \tilde{I}(\tilde g)>1+2\delta.
\end{equation}
Suppose that (\ref{el2.3}) is not true. Then there exists
$\{\tilde g_n, n\ge 1\}\subseteq \Theta_\rho$ such that $\lim_n
\tilde{I}(\tilde g_n)=1$. For $n$ large enough, $\tilde g_n$
belongs to the compact set $\{\tilde g;\ \tilde{I}(\tilde g)\le
2\}$, and there exists a subsequence $\{\tilde g_{n_k}, k\ge 1\}$
converging to $\tilde g$ in $\Theta_\rho$. As $\tilde{I}$ is lower
semicontinuous
$$1=\liminf_{k\rightarrow +\infty} \tilde{I}(\tilde g_{n_k})\ge \tilde{I}(\tilde
g),$$ and we get that $\tilde g \in \Theta$. So $d(\tilde g,
\Theta)=0$ what is a contradiction with the fact that $\tilde g
\in \Theta_\rho.$ Therefore, we can assume (\ref{el2.3}).

Now using (\ref{elimsup}) and (\ref{el2.3}) we have
\begin{equation*}
\limsup_{u\rightarrow + \infty} \frac{1}{\log\log u} \log
P\{\xi^u_\cdot\in \Theta_\rho\}\le - (1+2\delta),\end{equation*}
then, for $i$ large enough,
$$P\{\xi^{c^i}_\cdot\in \Theta_\rho\}\le
\exp\left\{-(1+2\delta)\log\log c^i\right\}\le
\frac{C}{i^{1+\delta}},$$ for some positive constant $C$. Finally,
the lemma is an immediate consequence of the Borel-Cantelli
lemma.\hfill $\Box$

\medskip

\noindent For every $i\ge 1$ and $c>1$ such that $c^{i-1}>e$,
define
$$\Gamma_i=\sup_{c^{i-1}\le u \le c^i}d \left(\xi^u_\cdot,
\frac{\phi(c^i)}{\phi(u)}\ \xi^{c^i}_\cdot \right)$$
%%\left(\frac{\eta^u_{\cdot}(X^u_0)}{\phi(u)},
%\frac{\eta^{y^i}_{\cdot}(X_0^{y^i})}{\phi(u)} \right).$$
\begin{lemma}\label{l2.4} For every $\rho>0$ there exists $c_\rho>1$
such that for $c\in (1,c_\rho)$ we have
$$P\left\{\exists \ i_0(\omega)\ {\rm s.t.}\ \Gamma_ i<\rho \ {\rm
whenever}\ i>i_0\right\}=1.$$
\end{lemma}
{\bf Proof}: We will prove that $$P(\limsup_{i \to \infty}
\{\Gamma_i\ge \rho\})=0.$$ From Lemma \ref{l2.3} it is enough to
check that
$$P(\limsup_{i \to \infty} \{\Gamma_i\ge \rho\}, \|\xi^{c^i}\|\le C)=0,$$ for
some positive constant $C$.

First,  notice that $$ \eta_t^u (X_0^u) = \frac{1}{\phi(u)}
\varphi_{ut}(X_0).$$ Then, for every $\delta>0$, taking $i$ large
enough, we get that $\frac{\phi(c^i)}{\phi(c^{i-1})} \le \sqrt{c
}(1+\delta)$ and so by means of this fact and using that $\phi$ is
nondecreasing we have, if $c$ is small enough,
\begin{equation}\label{edeixat}\begin{array}{l} \displaystyle\{\Gamma_i\ge
\rho, \|\xi^{c^i}\|\le C\}=\left\{ \sup_{c^{i-1}\le u \le c^i}d
\left(\frac{\varphi_{u \cdot}(X_0)}{\phi(u)}, \frac{\varphi_{c^i
\cdot}(X_0)}{\phi(u)}
\right)\ge \rho, \ \|\xi^{c^i}\|\le C\right\}\\[5mm]
\displaystyle\qquad \subseteq \left\{\sup_{0\le s\le 1}
\sup_{\frac{s}{c}\le t\le s} \frac{\phi(c^i)}{\phi(c^{i-1})}
\left| \frac{\varphi_{c^i t}(X_0)}{\phi(c^i)}- \frac{\varphi_{c^i
s}(X_0)}{\phi(c^i)}\right|\ge \rho, \ \|\xi^{c^i}\|\le
C\right\}\\[5mm]
\displaystyle\qquad \subseteq \left\{\sup_{0\le s\le 1}
\sup_{\frac{s}{c}\le t\le s} \left| \xi^{c^i}_t -
\xi^{c^i}_s\right|\ge \frac{\rho}{2}, \ \|\xi^{c^i}\|\le
C\right\}\\[5mm]
\displaystyle\qquad =\{\xi_\cdot^{c^i}\in
\Delta_\rho\},\end{array}\end{equation} where
$$\Delta_\rho=\left\{\tilde g \in
\mathcal{C}([0,1];\mathbb{R}^d); \sup_{0\le s\le 1}
\sup_{\frac{s}{c}\le t\le s} |\tilde g_t- \tilde g_s| \ge
\frac{\rho}{2}, \|\tilde g\|\le C\right\}.$$ Consider $f \in H^k$
such that $\tilde g=\tilde F_0(f)$ and $\tilde g\in \Delta_\rho$.
By (\ref{edetermnostre}), there exist $s\in [0,1]$ and $t\in
[\frac{s}{c},s)$ such that
\begin{equation*}
\sum_{j=1}^k \int_s^t |\tilde A_j(\tilde g_v)| | \dot f^j_v|\ dv
+\int_s^t  |\tilde A_0(\tilde g_v)| dv\ge | \tilde g_t- \tilde
g_s| \ge \frac{\rho}{4}.\end{equation*} On the other hand, using
the hypothesis of the coefficients (H) and assuming that $c<2$, we
have
$$\sum_{j=1}^k \int_s^t |\tilde A_j(\tilde g_v)|\ | \dot f^j_v|\ dv
+\int_s^t  |\tilde A_0(\tilde g_v)| dv \le C_1 \sqrt{2|t-s|\
\mathcal{I}(f)}+C_2 |t-s|,$$ for some positive constants $C_1$ and
$C_2$. Therefore
$$\mathcal{I}(f) \ge \frac{1}{C_1\sqrt{ 2(t-s)}} \left(\frac{\rho}{4}-C_2(t-s)\right),$$
 and this implies the existence of $c_\rho>1$ such
that, if $c\in (1,c_\rho)$, then $\mathcal{I}(f)>2$. Then
$$\inf\{I(\tilde g);\ \tilde g \in \Delta_\rho\}\ge 2.$$
Finally, for $i$ large enough, since $\Delta_\rho$ is closed, the
last estimate together with (\ref{edeixat}) and (\ref{elimsup})
yield that
$$P(\Gamma_i\ge \rho, \|\xi^{c^i}\|\le C)\le P(\xi_\cdot^{c^i}\in
\Delta_\rho)\le \exp\left\{-2\log\log c^i\right\}\le
\frac{C}{i^{1+\tau}},$$ for some $\tau>0$, and we can conclude
this lemma by means of the Borel-Cantelli lemma.\hfill $\Box$

\medskip

\begin{lemma}\label{l2.5} For every, $\rho>0$ there exists a.s. $u_0(\omega)>e$
such that, for every $u\in (u_0,+\infty)$, we have
$$d(\xi^u,\Theta)<\rho.$$
\end{lemma}
{\bf Proof}: Let $c>1$ and $i\in \mathbb{N}$ such that
$e<c^{i-1}<u\le c^i$, the triangular inequality gives
\begin{equation}\label{el2.5}
d(\xi^u,\Theta)\le d \left(\xi^u,\frac{\phi(c^i)}{\phi(u)}\
\xi^{c^i} \right)+d \left(\frac{\phi(c^i)}{\phi(u)}\
\xi^{c^i},\xi^{c^i}\right)+
d(\xi^{c^i},\Theta):=\beta_1+\beta_2+\beta_3.\end{equation} We
first deal with $\beta_1$.  Taking $c\in (1,+\infty)$  close to
$1$ and $i$ large enough, Lemma \ref{l2.4} yields
\begin{equation}\label{el2.5c} \beta_1=d \left(\xi^u,\frac{\phi(c^i)}{\phi(u)}\
\xi^{c^i} \right)<\frac{\rho}{3}.\end{equation} Study now
$\beta_2$. Lemma \ref{l2.3} implies that, for $i$ large enough,
$\| \xi_t^{c^i}\|$  is bounded. For every $\delta>0$ there exists
$i$ large enough such that $$1\le \frac{\phi(c^i)}{\phi(u)}\le
\frac{\phi(c^i)}{\phi(c^{i-1})}\le \sqrt{c}(1+\delta).$$ Then,
using these two facts, for $c$ close enough to 1, we have that
\begin{eqnarray}
\beta_2&=&d\left(\frac{\phi(c^i)}{\phi(u)}\
\xi^{c^i},\xi^{c^i}\right)= \sup_{0\le t\le 1} \left|
\frac{\phi(c^i)}{\phi(u)}\
\xi^{c^i}_t-\xi^{c^i}_t\right|\nonumber\\
&\le &\sup_{0\le t\le 1} |\xi^{c^i}_t|\
\left|\frac{\phi(c^i)}{\phi(u)}\
-1\right|<\frac{\rho}{3}.\label{el2.5b}
\end{eqnarray} Lemma \ref{l2.3} implies that, for $i$ large
enough,
\begin{equation}\label{el2.5a}
\beta_3=d(\xi^{c^i},\Theta)<\frac{\rho}{3}.\end{equation} So, we
finish the proof of this lemma applying
(\ref{el2.5c})-(\ref{el2.5a}) to  (\ref{el2.5}). \hfill $\Box$

\medskip

\begin{lemma}\label{l2.6} Consider $\tilde g \in \Theta$ such that
$\tilde{I}(\tilde g)<1$. Then, for every $\rho>0$, there exists
$c_\rho>1$ such that, for every $c>c_\rho$, we have
$$P\left\{d(\xi^{c^i},\tilde g)<\rho, {\rm infinitely\,
often}\right\}=1.$$
\end{lemma}
{\bf Proof}:  Let $f\in H^k$ such that $\tilde F_0(f)=\tilde g$
and $\mathcal{I}(f)<1$. For fixed $\rho>0$, define for $\nu>0$
$$\Upsilon_i=\left\{\left\|\frac{1}{\sqrt{c^i \log\log c^i}}\
W_{c^i\cdot}-f\right\|\le \nu\right\}\qquad {\rm and} \qquad
\Lambda_i=\big\{\big\|\xi^{c^i}-\tilde g\big\| \le \rho\big\}.$$
Since $P(\limsup_{i \to \infty} \Upsilon_i)=1$, following the same
argument as in  Lemma 2.6 of Baldi (1986) and as a consequence of
the scaling property, we only need to prove that
\begin{equation}\label{sumafinita}
\sum_i P(\Upsilon_i\cap \Lambda_i^c) < \infty.
\end{equation}
Notice first that
\begin{equation*}
P(\Upsilon_i\cap \Lambda_i^c) =P\left(\left\| \xi^{c^i} - \tilde g
\right\|>\rho, \left\|\frac{1}{\sqrt{\log\log c^i}}W - f \right\|
\le \nu\right)\le P_1+P_2,\end{equation*} with
\begin{eqnarray*}
P_1&=&  P\left(\left\| \xi^{c^i} - \tilde g \right\|
>\rho, \left\|\frac{1}{\sqrt{\log\log
c^i}}W - f \right\| \le \nu,   \vert X_0^{c^i} \vert \le  \tau \right),\\
P_2&=&P\left(\vert X_0^{c^i} \vert \ge  \tau\right).
\end{eqnarray*}
For any $\tau$, set $O_\tau$ the closed ball $B(0,\tau)$. Fixed
$\rho$ and $\tau$, using Theorem \ref{desnostre}, there exits $\nu
>0$ and $i_0$ such that, for any $i \ge  i_0$,
\begin{eqnarray}
P_1 & \le & P\left(   \left\| \eta^{c^i} - \tilde g
\right\|_{O\tau}
>\rho, \left\|\frac{1}{\sqrt{\log\log
c^i}}W - f \right\| \le \nu \right) \nonumber
\\
& \le & \exp (-2 \log \log c^i ) \le \frac{C}{i^2}.\label{desp1}
\end{eqnarray}
On the other hand, hypothesis (C) yields
\begin{equation}
P_2 \le \exp (-2 \sqrt{\log \log c^i} ) \le
\frac{C}{i^2},\label{desp2}
\end{equation}
for $i$ big enough.

Putting together (\ref{desp1}) and (\ref{desp2}), we easily obtain
(\ref{sumafinita}). \hfill $\Box$

\bigskip

\noindent {\bf Proof of Theorem \ref{tlil}}: Lemma \ref{l2.5}
implies that $\{ \xi^u \}_u$ is relatively compact. Moreover,
Lemma \ref{l2.6} ensures that all the points of $\Theta$ are limit
points.

\hfill $\Box$

\end{document}